\documentclass[12pt]{article} 
{  
\begin{document} 
%%%%%%%%%%%%%%%%%%%%%%
%%% ** start of amsfont definitions **
\newfont{\twelvemsb}{msbm10 scaled\magstep1}
\newfont{\eightmsb}{msbm8} \newfont{\sixmsb}{msbm6} \newfam\msbfam
\textfont\msbfam=\twelvemsb \scriptfont\msbfam=\eightmsb
\scriptscriptfont\msbfam=\sixmsb \catcode`\@=11
\def\Bbb{\ifmmode\let\next\Bbb@\else \def\next{\errmessage{Use
\string\Bbb\space only in math mode}}\fi\next}
\def\Bbb@#1{{\Bbb@@{#1}}} \def\Bbb@@#1{\fam\msbfam#1}
\newfont{\twelvegoth}{eufm10 scaled\magstep1}
\newfont{\tengoth}{eufm10} \newfont{\eightgoth}{eufm8}
\newfont{\sixgoth}{eufm6} \newfam\gothfam
\textfont\gothfam=\twelvegoth \scriptfont\gothfam=\eightgoth
\scriptscriptfont\gothfam=\sixgoth \def\frak{\frak@}
\def\frak@#1{{\fam\gothfam{{#1}}}} \def\frak@@#1{\fam\gothfam#1}
\catcode`@=12
%%% ** end of amsfont definitions **
%
%\def\Bbb{\bf}
\def\CC{{\Bbb C}}
\def\NN{{\Bbb N}}
\def\QQ{{\Bbb Q}}
\def\RR{{\Bbb R}}
\def\ZZ{{\Bbb Z}}
\def\cA{{\cal A}} \def\cB{{\cal B}} \def\cC{{\cal C}}
\def\cD{{\cal D}} \def\cE{{\cal E}} \def\cF{{\cal F}}
\def\cG{{\cal G}} \def\cH{{\cal H}} \def\cI{{\cal I}}
\def\cJ{{\cal J}} \def\cK{{\cal K}} \def\cL{{\cal L}} 
\def\cM{{\cal M}} \def\cN{{\cal N}} \def\cO{{\cal O}}
\def\cP{{\cal P}} \def\cQ{{\cal Q}} \def\cR{{\cal R}} 
\def\cS{{\cal S}} \def\cT{{\cal T}} \def\cU{{\cal U}}
\def\cV{{\cal V}} \def\cW{{\cal W}} \def\cX{{\cal X}}
\def\cY{{\cal Y}} \def\cZ{{\cal Z}}
\def\qed{\hfill \rule{5pt}{5pt}}
\def\th{\mbox{\footnotesize th}}
\newenvironment{result}{\vspace{.2cm} \em}{\vspace{.2cm}}
%
%%%%%%%%%%%%%
\begin{center} 
{\mathversion{bold}
{\LARGE {\bf A two-parametric deformation\\[2mm] of $U[sl(2)]$, its representations and complex "spin"}}}
\\[9mm] 
{\bf Nguyen Anh Ky} and {\bf Nguyen Thi Hong Van}\\[2mm] 
Institute of Physics and Electronics\\ 
P.O. Box 429, Bo Ho, Hanoi 10 000, Vietnam 
\end{center} 
\vspace*{1mm} 
\begin{abstract}     
A two-parametric deformation of  $U[sl(2)]$ and its representations  are considered. This 
newly introduced two-parametric quantum group  denoted as $U_{p,q}[sl(2)]$ admits a 
class of infinite - dimensional representations which have no classical (non-deformed) 
and one-parametric deformation analogues, even at generic deformation parameters. 
Interestingly that finite - dimensional representations of $U_{p,q}[sl(2)]$ allow arbitrary 
complex "spins" (i.e., not necessary they to be integral or half-integral numbers), unlike 
those in the classical and one-parametric deformation cases. 
\end{abstract} 
\begin{center}
\end{center}
\vspace*{5mm}
{\Large {\bf 1. Introduction}} 
 
Quantum groups  are one of the most fascinating mathematical concepts with a 
physical origin \cite{frt} - \cite{woro}. 
Depending on points of view,  quantum groups  can be approached to in several ways. 
One of the approaches to quantum groups are the so-called Drinfel'd-Jimbo deformation 
of universal enveloping algebras \cite{drin, jim}. The quantum groups  of this kind 
(also called quantum  algebras) are non-commutative and non-cocommutative Hopf 
algebras \cite{drin}. 
By construction, such a quantum group depends on a, complex in general, 
parameter called a (quantum) deformation parameter.  It is logical to ask if these 
one-parametric deformations can be extended to multi-parametric ones. This question 
was considered first by Manin \cite{manin1} and later by a number of authors 
(see \cite{sudbe} - \cite{ck} and references therein).  The structure of a multi-parametric 
deformation is usually richer than that of an one-parametric deformation. Unfortunately, 
in comparison with one-parametric deformations, multi-parametric deformations are less 
understood, in spite of some progress made in this direction. Therefore, the latter deserve 
more comprehensive investigations in both the mathematical and physical aspects. 
Here we consider  a relatively simple case, a two-parametric deformation of $U[sl(2)]$, 
which, however, has interesting structure and features. 

The quantum  group  $U_q[sl(2)]$ as an one-parametric deformation of the universal 
enveloping algebra $U[sl(2)]$ is one of the best  investigated quantum groups 
\cite{kure} - \cite{pas}. What about two-parametric deformations of $U[sl(2)]$, they 
have been considered  in several versions and in different aspects (see, for example, 
\cite{swz} - \cite{pq96} and references therein) though some of them  are, in fact, 
equivalent to one-parametric deformations (upto some rescales).  In this paper we 
introduce and consider one more two-parametric deformation of $U[sl(2)]$ denoted 
as $U_{p,q}[sl(2)]$. Following the method of highest-weight  representations we can 
contruct its representations. It turns out that this new quantum group admits  a class 
of infinite-dimensional representations which have no analogue in the cases of the 
non-deformed $sl(2)$ and previously introduced deformations of $U[sl(2)]$. 
Morever, the "spin" (highest weight) corresponding to a finite-dimensional 
representation of $U_{p,q}[sl(2)]$ could be an arbitrary complex number, unlike 
the finite-dimensional representations of $sl(2)$  and its one-parametric deformations 
for which a "spin" is an (half) integral number. We note that $U_{p,q}[sl(2)]$ is not 
equivalent to the one-parametric $U_q[sl(2)]$ unless at very special choices of $p$ 
and $q$ such as $p=q$.
 
In the next section, for compare, we briefly recall the non-deformed $sl(2)$ and 
its representations. The quantum group $U_{p,q}[sl(2)]$ is introduced and considered 
in Sect. 3.  Some discussions and conclusions are made in the last section, Sect.4. \\[7mm]
{\Large {\bf 2. {\mathversion{bold}$sl(2)$} and representations}} 

The algebra $sl(2)$ can be generated by three generators, say $E_+$, $E_-$ and $H$,  subject to  
the commutation relations 
$$[H,E_{\pm}]=\pm E_{\pm},  \qquad [E_+,E_-]=2H.\eqno(1)$$ 
Demanding $(H)^{\dagger}=H$ and $(E_{\pm})^{\dagger}=E_{\mp}$, a (unitary) 
representtaion induced from  a (normalised) highest weight state $|j,j\rangle$ with 
a highest weight ("spin") $j$, 
$$H|j,j\rangle=j|j,j\rangle, \quad E_+|j,j\rangle=0,\eqno(2)$$ 
has the matrix elements 
$$H~|j,m\rangle = m ~|j,m\rangle, $$ 
$$ E_\pm~| j,m\rangle = 
\sqrt{(j\mp m)(j\pm m +1)}~ |j, m\pm 1\rangle,\eqno(3)$$ 
where $|j,m\rangle$, $m \leq j$, is one of the orthonormalised states, 
$$\langle j,m_1|j,m_2\rangle =\delta_{m_1m_2}, $$
obtained from the highest weight state $|j,j\rangle$ by acting on the latter  a monomial 
of the generator $E_-$ of an appropriate order, say $n$, 
$$|j,m\rangle= A_n(E_-)^n~|j,j\rangle,\qquad n\in \NN, \quad m=j-n, \eqno(4)$$ 
with $A_n$ a normalizing coefficient, which for a finite-dimensional  representation 
(i.e., for a non-negative (half) integral $j$) equals
$$A_n=\sqrt{(2j)!\  n!\over (2j-n)!}\  .\eqno(5)$$ 
 In this case, the representations constructed  are simultaneously highest weight  
(above-bounded) and lowest weight (below-bounded), that is, 
 $$E_+|j,j\rangle=0, \quad  E_-|j,-j\rangle=0.\eqno(6)$$ 
They are finite-dimensional (and also unitary and irreducible) representations of 
dimension $2j+1$. The situation is similar in the case of the one-parametric quantum 
group $U_q[sl(2)]$ at a generic deformation parameter $q$ (i.e., at $q$ not a root of 
unity)  
\cite{ss03}. \\[7mm] 
{\Large {\bf 3. Two-parametric deformation {\mathversion{bold}$U_{p,q}[sl(2)]$}}} 

The two-parametric quantum group $U_{p,q}[sl(2)]$ as a two-parametric deformation of 
$U[sl(2)]$ is generated also 
by three generators  $E_+$, $E_-$ and $H$ which now satisfy the deformed defining 
relations 
$$[H,E_{\pm}]=\pm E_{\pm},  \quad [E_+,E_-]=[2H]_{p,q},\eqno(7)$$
where $$[x]_{p,q}=\frac{q^x-p^{-x}}{q-p^{-1}}\eqno(8)$$ 
is a two-parametric deformation of $x$ (a number or an operator) with $p$ and $q$ 
being complex, in general, deformation parameters ($p^2 \neq q^2$). When $p=q$ 
this two-parametric deformation is reduced to the one-parametric deformation 
$U_{q}[sl(2)]$.
 
At generic $p$ and $q$, the center of $U_{p,q}[sl(2)]$ is spanned on the Casimir operator 
$$C=\frac{1}{1-q^{-2}}(k_q)^2-\frac{1}{1-p^2}(k_p^{-1})^2+(q-p^{-1})E_-E_+,
\eqno(9)$$
where $k_q:=q^H$, $k_p^{-1}:=p^{-H}$. The latter can be used to quickly construct unitary 
representations of $U_{p,q}[sl(2)]$. We find the representations of $U_{p,q}[sl(2)]$ 
corresponding to those of  $sl(2)$ in (3),   
$$H~|j,m\rangle = m ~|j,m\rangle, $$ 
$$E_+~|j,m\rangle = 
\left ( ~ {q^{j+m+1} [j-m]_q-p^{-j-m-1} [j-m]_p}\over {q-p^{-1}}\right ) ^{1/2} 
~ |j, m+1\rangle,$$ 
$$\quad E_-~|j,m\rangle =  
\left  ( ~ {q^{j+m}[j-m+1]_q-p^{-j-m}[j-m+1]_p}\over {q-p^{-1}}\right ) ^{1/2} 
~ |j, m-1\rangle.\eqno(10)$$ \\ 
When the unitary condition $(E_{\pm})^{\dagger}=E_{\mp}$ is imposed, the coefficient $A_n$ 
satisfies, instead of (4), the recurrent formula
$$\left |{A_{n-1}}\over {A_n}\right |^2=
\left. {q^{2j-n+1}[n]_q-p^{-2j+n-1}[n]_p}\over {q-p^{-1}}\right., ~~~ n=j-m$$
if the latter is meaningful (i.e., if its r.h.s is a real positive number for some choice 
of $j$, $p$ and $q$). Generally speaking, however, 
the representations constructed are not unitary and we do not need this recurrent 
formula to obtain (10). These representations of 
$U_{p,q}[sl(2)]$, even at  (half) integral $j$,  are, in general, infinite-dimensional, 
unlike the constructed in a similar way representations of $sl(2)$ and other its 
deformations \cite{ss03} which are finite-dimensional for non-negative (half) integral 
$j$. At arbitrary $p$ and $q$, the representaions (10) of $U_{p,q}[sl(2)]$ are highest 
weight (by construction) but, as can be seen from (10),  they may not be lowest weight 
(unbounded from below), i.e., not finite-dimensional 
any more, even for an integral $2j$. Therefore, it is a new class of infinite-dimensional 
representations of $U_{p,q}[sl(2)]$ not found before in the cases of $sl(2)$ 
and its previously considered deformations. \\[4mm]
{\bf Proposition 1}: {\it Highest weight representations of the two-parametric quantum 
group $U_{p,q}[sl(2)]$ given in (10) are in general infinite-dimensional, even for 
non-negative (half) integral highest weights}.

The next interesting phenomenon is related to the finite-dimensional representations. 
The representations (10) are finite-dimensional if the matrix element of $E_-$ vanishes 
for some value of $n\equiv j-m$ which is a non-negative integer ($n\in \NN$). \\[4mm]
{\bf Proposition 2}: {\it  For a given, not necessary (half) intergral, $j$, the 
representation (10) is reducible and contains a finite-dimensional subrepresentation 
iff the equation 
$$f(x)\equiv  q^{2j-x}[x+1]_q-p^{-2j+x}[x+1]_p=0\eqno(11)$$
has a non-negative integral solution.}\\[4mm]
The request for an integral solution makes the last equation resembling an equation of 
Diophantine type (the difference here is the coefficients in (11) are not necessarily 
integral).  In general, to prove the Eq. (11) to have or not an integral solution (and when) 
is a hard mathematical problem  which needs futher investigations.  Suppose $x=\cN$ is 
the smallest non-negative integer solving Eq. (11), the dimension of the corresponding 
finite-dimensional reprsentation extracted from (10) is $D=\cN+1$ which may not 
equal $2j+1$. Contrarily, we can choose $j$ to get a representation of a given finite 
dimension. \\[4mm]
{\bf Proposition 3}: {\it  For a representation of a given finite dimension $D$ the highest 
weight is given by} 
$$2j=\log_{pq}\left (\frac{[D]_p}{[D]_q}\right)+D-1 \quad \mbox{or} \quad 
2j=\frac{\ln\left (\frac{[D]_p}{[D]_q}\right)}{\ln (pq)}+D-1\eqno(12)$$
{\it if the logarithms here are well-defined.}
\\[4mm]
When the logarithms in (12) are  not well-defined we should keep the Eq. (11) and solve 
it for $j$ with a given value of $x=D-1\in \NN$. The feature is $2j$, if found, may not be a 
non-negative integer at all, but an arbitrary complex number, $2j\in \CC$. 

The representations (10) have a relatively simple structure at generic $p$ and $q$. 
In the case of one or both parameters being roots of unity, the Casimir operator (9) 
may not exhaust an (extended) center of $U_{p,q}[sl(2)]$ anymore and the representation 
structure becomes in general more complex (but sometimes the representation structure is 
similar to that at generic parameters). This interesting case deserves to be separately 
investigated in details.   

Note that $U_{p,q}[sl(2)]$ is by no means equivalent to the one-parametric $U_q[sl(2)]$  
(unless at very special choices of $p$ and $q$) and other previously deformations of 
$U[sl(2)]$. \\[5mm] 
{\Large 4. {\bf Conclusion}} 

The quantum group $U_{p,q}[sl(2)]$ which is a two-parametric deformation of $U[sl(2)]$ 
has been introduced and its representations have been investigated. This quantum group, 
even at generic deformation parameters, has intersting features. It it is  showed that 
$U_{p,q}[sl(2)]$  admits a class of infinite-dimensional representations which have no 
analogues in the case of $sl(2)$ and its prviously considered deformations \cite{ss03}.  
It is a new phenomenon of $U_{p,q}[sl(2)]$. In gereral, the representations found are  
irreducible. They may become reducible under certain circumstances and then 
finite-dimensional irreducible subrepresentations can be extracted.  Another feature 
is the highest weight (the "spin") of a finite-dimensional representation of 
$U_{p,q}[sl(2)]$ is not necessarily a (half) integral number but a complex one, 
unlike that of $sl(2)$ and its one-parametric deformation. This fact may have an 
interesting physical interpretation. Other classes of infinite - dimensional 
representations of $U_{p,q}[sl(2)]$ and applications may be found by using the 
methods of \cite{stoyanov,kundu}. 

The next problem one could consider is $U_{p,q}[sl(2)]$ at roots of uinity. It is 
also interesting to investigate this quantum group and its co-algebra structure 
in the light of the $R$-matrix formalism and to look for possible associated 
integrable models. Some of these problems turn out to be harder than expected  
at first sight. In this spirit, the results obtained in the present paper are 
far from being complete but we think they could be of independent interest.  
\\[4mm]  
{\bf Acknowledgement}: This work was partially supported by the Vietnam National 
Research Program for Natural Sciences under Grant No  410804.  
 
\end{document}